\documentstyle[amsmath,amssymb,twoside,12pt]{article}
\allowdisplaybreaks

\def\Sum{\displaystyle\sum}
\def\Frac{\displaystyle\frac}

\def\R{{\mathbb R}}

\begin{document}
\begin{center}
{\Large\bf A new method for the computation of eigenvalues}
\vskip 1cm {\bf Nassim  Guerraiche }\\
\
 Laboratoire des Math\'ematiques Appliqu\'es et Pures, Universit\'e de
Mostaganem\\ B.P. 227, 27000, Mostaganem, ALG\'ERIE\\
nassim.guerraiche@univ-mosta.dz \\[0.3cm]
\vskip0.2cm
\end{center}

\begin{abstract}

In this paper we are concerned to find the  eigenvalues and eigenvectors of a real symetric matrix  by applying a new numerical method similar to Jacobi method. Our approch consists to use a new  orthogonal matrix. The computation of the eigenvalues and eigenvectors by using this method appears easier if compared with Jacobi method in the sense of the functions used in the orthogonal matrix.
\end{abstract}



\section{Introduction}
As we know, for a given  matrix $A\in \mathbb{C}^{\text{n} \times \text{m}}$,  the computation of its eigenvalues and eigenvectors  is easy when its dimension is smal but this computation will become difficult when the dimension of the matrix is big. For the  matrices with big dimensions many researchers have contributed and gave  different numerical methods to compute their eigenvalues and eigenvectors, as exapmles,  we cite the QR method (see \cite{kubla, fran}), the power method \cite{muntz} and Sturm sequences method which can be found in the book of Quarteroni {\em et al} \cite{quart}; and for a real symetric matrix usually we use the Jacobi method \cite{jacobi}.\\

This paper deals with the computation of  eigenvalues and eigenvectors of a real symetric matrix $A\in \mathbb{R}^{\text{n} \times \text{n}}$, by changing the Givens matrix using on the Jacobi method by another orthogonal matrix, i.e.,  we will replace the matrix: 
\[ G = \begin{pmatrix}
1 & 0 &  &  & &0\\
0 & \ddots &  & &  & \\
 &  & \cos(\theta) & -\sin(\theta)& & \\
 &  & \sin(\theta) & \cos(\theta)& & \\
& & & &\ddots&\\
0& & & & &1
\end{pmatrix} \]
by the following matrix H:  
\[ H = \begin{pmatrix}
1 & 0 &  &  & &0\\
0 &  &  & &  & \\
 &  & \sqrt{x+\delta} & -\sqrt{-x-\delta+1}& & \\
 &  & \sqrt{-x-\delta+1} &\sqrt{x+\delta} & & \\
& & & & &\\
0& & & & &1
\end{pmatrix} \]
such that \quad $\delta \in \R$  \quad and \quad  $-\delta \leq x \leq 1-\delta$.

\vspace{0.5cm}

The paper is organized as follows: in Section $2$ we present our main results and in Section $3$ we give a MATLAB program to this new method.  In the remaider of this paper and without loss of generality  we will choose $\delta=\frac{1}{2}$.\\

\section{The computation of eigenvalues and eigenvectors}

\vspace{0.3cm}
In this section, in order to find the eigenvalues and eigenvectors of a real symetric matrix $A\in \mathbb{R}^{\text{n} \times \text{n}}$,  we will repeat the procedures of Jacobi method by using the matrix $H$ introduced above and we shall give all the steps of the calculus. \\

Let $A\in \mathbb{R}^{\text{n} \times \text{n}}$ be a real symetric matrix of $n$ dimensions with coefficients $(a_{ij})_{1 \leq i,j \leq n}$.\\
The jacobi method \cite{quart} is an iterative method consists to build a sequences of orthogonal matrices $A^{(k)}$ such that on the $k$-iteration we have
\begin{center} 
$B^{(k)}=H_{pq}^{T} A^{(k-1)} H_{pq},\quad (A^{(0)}=A)$
\end{center}
where $a^{(k)}_{ij}=0$  if  $(i,j)=(p,q)$ and the matrix $B^{(k)}$ converge to the matrix of eigenvalues. \\

\vspace{0.2cm}

In a block $A^{(k-1)}$ of the matrix $A$, we find
\[
\mathbf{\begin{pmatrix}
a^{(k)}_{pp} & a^{(k)}_{pq} \\  a^{(k)}_{pq}& a^{(k)}_{qq} \end{pmatrix}} = \begin{pmatrix}
\sqrt{x+\frac{1}{2}} & \sqrt{-x+\frac{1}{2}} \\ -\sqrt{-x+\frac{1}{2}} & \sqrt{x+\frac{1}{2}} \end{pmatrix}
\begin{pmatrix}
a^{(k-1)}_{pp} & a^{(k-1)}_{pq} \\  a^{(k-1)}_{pq}& a^{(k-1)}_{qq} \end{pmatrix}
\begin{pmatrix}
\sqrt{x+\frac{1}{2}} & -\sqrt{-x+\frac{1}{2}} \\\sqrt{-x+\frac{1}{2}}  & \sqrt{x+\frac{1}{2}} \end{pmatrix}
\]
\[
= \begin{pmatrix}
a_{pp}(x+\frac{1}{2})+2a_{pq}\sqrt{-x^{2}+\frac{1}{4}}+(-x+\frac{1}{2})a_{qq}* & (-a_{pp} +a_{qq})\sqrt{-x^{2}+\frac{1}{4}}+2a_{pq}x \\ (-a_{pp}+a_{qq})\sqrt{-x^{2}+\frac{1}{4}}+2a_{pq}x & a_{pp}(-x+\frac{1}{2})+2a_{pq}\sqrt{-x^{2}+\frac{1}{4}}+(x+\frac{1}{2})a_{qq}** \end{pmatrix}
\]

where $1 \leq p < q \leq n$\\

Now, when we solve the equation

\begin{equation}\label{e1}
(a_{qq}-a_{pp})\sqrt{-x^{2}+\frac{1}{4}}+2a_{pq}x=0
\end{equation} 

we find that

\begin{equation*}
x_{0}=\pm \frac{|a_{qq}-a_{pp}|}{2\sqrt{(a_{qq}-a_{pp})^{2}+4a^{2}_{pq}}}
\end{equation*}

i.e.,  if $a_{pq} >0$, we have

\begin{equation*}
x_{0}=\frac{a_{pp}-a_{qq}}{2\sqrt{(a_{qq}-a_{pp})^{2}+4a^{2}_{pq}}}
\end{equation*}

then by subtitution in $(*)$ and $(**)$ we find

\begin{equation*}
\lambda_{*}=\frac{a_{qq}+a_{pp}}{2}+\frac{{(a_{pp}-a_{qq})}^{2}+4a^{2}_{pq}}{2\sqrt{(a_{qq}-a_{pp})^{2}+4a^{2}_{pq}}}
\end{equation*}

\begin{equation*}
\lambda_{**}=\frac{a_{qq}+a_{pp}}{2}+\frac{-{(a_{pp}-a_{qq})}^{2}-4a^{2}_{pq}}{2\sqrt{(a_{qq}-a_{pp})^{2}+4a^{2}_{pq}}}
\end{equation*}

and if  $a_{pq} <0$, the root of the equation is 

\begin{equation*}
x_{0}=\frac{a_{qq}-a_{pp}}{2\sqrt{(a_{qq}-a_{pp})^{2}+4a^{2}_{pq}}}
\end{equation*}

by subtitution $x_{0}$ by its value in $(*)$ and $(**)$ we find

\begin{equation*}
\lambda_{*}=\frac{a_{qq}+a_{pp}}{2}+\frac{-{(a_{pp}-a_{qq})}^{2}-4a^{2}_{pq}}{2\sqrt{(a_{qq}-a_{pp})^{2}+4a^{2}_{pq}}}
\end{equation*}
\begin{equation*}
\lambda_{**}=\frac{a_{qq}+a_{pp}}{2}+\frac{{(a_{pp}-a_{qq})}^{2}+4a^{2}_{pq}}{2\sqrt{(a_{qq}-a_{pp})^{2}+4a^{2}_{pq}}}
\end{equation*}

\subsection{Existence and uniqueness of  solution of the equation \ref{e1}}

In this section we show that the equation (\ref{e1}) has a unique solution $x_{0}$. The idea of the proof consists to divide the  interval  $[-\frac{1}{2},\frac{1}{2}]$ on two open intervals,  $]-\frac{1}{2},0[$ and  $]0,\frac{1}{2}[$.\\ 

\vspace{0.5cm}

Let 
\begin{center}
$f(x)=(a_{qq}-a_{pp})\sqrt{-x^{2}+\frac{1}{4}}+2a_{pq}x$
\end{center}
with the dervative
\begin{center}
$f^{\prime}(x)=(a_{qq}-a_{pp})\frac{-x}{\sqrt{-x^{2}+\frac{1}{4}}}+2a_{pq}$
\end{center}
now,  for $a_{qq}-a_{pp} >0$, $a_{pq}>0$ and  for $x \in ]-\frac{1}{2},0[$, we find

\begin{equation}\label{e2}
-a_{pq} < f(x) < \frac {a_{qq}-a_{pp}}{2}
\end{equation}
also we can find when $x \in ]0,\frac{1}{2}[$ the following

\begin{equation}\label{e3}
0 < f(x) < \frac {a_{qq}-a_{pp}}{2}+a_{pq}
\end{equation}
By (\ref{e2}) and according to the intermediate value theorem, $f$ has at least one root in the interval $ ]-\frac{1}{2},0[$ and by using the fact that $f^{\prime}$ is stricly  positive in this interval, we can say that the root is unique. Now the expression (\ref{e3}) shows that $f$ keeps its   positive sign on the interval $]0,\frac{1}{2}[$, it means that  $f$ do not have any roots in this interval, so  from the above we deduce that $f$ has a unique root in the interval $[-\frac{1}{2},\frac{1}{2}]$. 

\vspace{0.4cm}

Processing by the same manier we can find:
\begin{enumerate}
\item if $a_{pq} >0$ and $a_{qq}-a_{pp}<0$ we have respectively on the intervals $ ]-\frac{1}{2},0[$ and $]0,\frac{1}{2}[$
\begin{equation*}
\frac {a_{qq}-a_{pp}}{2}-a_{pq} < f(x) < 0
\end{equation*}
\begin{equation*}
\frac {a_{qq}-a_{pp}}{2} < f(x) < a_{pq}
\end{equation*}
 then $x_{0} \in ]0,\frac{1}{2}[$
\item if $a_{pq} <0$ and $a_{qq}-a_{pp}>0$ we have respectively on the intervals $ ]-\frac{1}{2},0[$ and $]0,\frac{1}{2}[$
\begin{equation*}
0 < f(x) < \frac {a_{qq}-a_{pp}}{2}-a_{pq}
\end{equation*}
\begin{equation*}
a_{pq} < f(x) < \frac {a_{qq}-a_{pp}}{2} 
\end{equation*}
 then $x_{0} \in ]0,\frac{1}{2}[$
\item if $a_{pq} <0$ and $a_{qq}-a_{pp}<0$  we have respectively on the intervals $ ]-\frac{1}{2},0[$ and $]0,\frac{1}{2}[$
\begin{equation*}
\frac {a_{qq}-a_{pp}}{2}< f(x) < -a_{pq}
\end{equation*}
\begin{equation*}
\frac {a_{qq}-a_{pp}}{2}+a_{pq} < f(x) < 0
\end{equation*}
then $x_{0} \in ]-\frac{1}{2},0[$
\end{enumerate}
\section{MATLAB Program }
Although this method is very similar to the Jacobi method, which is of course convergent, this does not prevent us  to give it an associate program. In what follows we shall give the program of this new method. Our approch is based directly upon the program of cyclic Jacobi method given in \cite{quart} (Program 23-33, 35-37). A few changes were made since the functions of the orthogonal matrix were changed.\\

It is clear that the numerical esimations of the Jacobi method is still here unchaged. First of all, Let give the following quantity
$$
\Psi(A)=\left(\Sum \limits_{\underset{i \not=j}{i, j=1}}^{n}a_{ij}^2\right)^{1/2}=\left(\|A\|_{F}^{2}-\Sum_{i=1}^{n}a_{ii}^{2}\right)^{1/2}
$$
such that $\|\cdot\|_{F}$ is the Frobenius norm. And it is well known that in the k-iteration we have
$$
\Psi(A^{(k)}) \leq \Psi (A^{(k-1)}), \, \, \, \text{  for  } k \geq 1 
$$
Let also give the following estimation
$$
\Psi(A^{(k+n)}) \leq \Frac{1}{\delta \sqrt{2}}(\Psi(A^{(k)}))^{2}, \, \, \, k=1, 2,\cdots
$$
this last is obteined in the cyclic Jacobi method, where $N=n(n-1)/2$ and $\delta$,  by hypothes,  satisfies the following inequality
$$
|\lambda_{i}-\lambda_{j}|\geq \delta \, \, \, \text{ for } i \not=j 
$$
Now, we give the MATLAB program with the changes required.
\begin{itemize}
\item Let start by the program that allows us to calculate the product $H(i,k,x)M$\\

function [M]=pro1(M,irr1,irr2,i,k,j1,j2)\\
for j=j1:j2\\
    t1=M(i,j);\\
    t2=M(k,j);\\
    M(i,j)=irr1.*t1+irr2.*t2;\\
    M(k,j)=-irr2.*t1+irr1.*t2;\\
end\\
return\\

such that $irr1=\sqrt{x+1/2}$ and $irr2=\sqrt{-x+1/2}$ 

\item Secondly, we give the program of the product $MH(i,k,x)^{T}$\\

function[M]=pro2(M,irr1,irr2,j1,j2,i,k)\\
for j=j1:j2\\
    t1=M(j,i);\\
    t2=M(j,k);\\
    M(j,i)=irr1*t1+irr2*t2;\\
    M(j,k)=-irr2*t1+irr1*t2;\\
end\\
return\\

\item Now, we give the program which allows us to  evaluate  $\Psi(A)$ in the cyclic new  method\\

function[psi]=psinorm(A)\\
$[n,m]$=\text{size}(A);\\
if n$\not=$m, error('only for square matrix'); end\\
psi$=$0;\\
for i=1:n-1\\
    j=$[i+1:n]$;\\
    psi=psi+sum($A(i,j).^2$+$A(j,i).^2$')\\
end\\
psi=sqrt(psi);\\
return\\

\item Afterwards, the program which allows us to  evaluate  $irr1$ and $irr2$\\

function[irr1,irr2]=symschur2(A,p,q)\\
if A(p,q)==0\\
    irr1=1;irr2=0;\\
else\\
    if A(p,q)$>=$0\\
        z1=(A(p,p)-A(q,q));\\
        z2=((A(q,q)-$A(p,p)).^2$)+(4.*$(A(p,q)).^2$);\\
        z3=sqrt(z2);\\
        z4=2.*z3;\\
        x=z1./z4;\\
    else\\
        v1=(A(q,q)-A(p,p));\\
        v2=((A(q,q)-$A(p,p)).^2$)+(4.*$(A(p,q)).^2$);\\
        v3=sqrt(v2);\\
        v4=2.*v3;\\
        x=v1./v4;\\
    end\\
    irr1=sqrt(x+(1/2)); irr2=sqrt(-x+(1/2));\\
end\\
return\\

\item Finally, here is the program of the new method\\

function[D,sweep,psi]=cycjacobi2(A,tol,nmax)\\
$[n,m]$=size(A);\\
if n$\not=$m, error('only for the square matrix'); end\\
D=A;\\
psi=norm(A,'fro');\\
epsi=tol*psi;\\
psi=psinorm(D);\\
sweep=0;\\
iter=0;\\
while psi$>$epsi and iter$<=$nmax\\
    iter=iter+1;\\
    sweep=sweep+1;\\
    for p=1:n-1\\
        for q=p+1:n\\
            $[irr1,irr2]$=symschur2(D,p,q);\\
            $[D]$=pro1(D,irr1,irr2,p,q,1,n);\\
            $[D]$=pro2(D,irr1,irr2,1,n,p,q);\\
        end\\
    end\\
    psi=psinorm(D);\\
end\\
return\\

such that tol is the tolerance and nmax is the maximum number of iterations.
\end{itemize}
\section{Example}
Let 
\[ A = \begin{pmatrix}
1 & 0 & 2&\\
0 & 3 &0&\\
 2&0  & 4 
\end{pmatrix} \]
and let $A^{(0)}=A$, then we have
\begin{center}
$A^{(1)}=H^{T}AH$
\end{center}
such that
\[ H = \begin{pmatrix}
\sqrt{x+\frac{1}{2}} &0& -\sqrt{-x+\frac{1}{2}}&\\
 0&1&0& \\
\sqrt{-x+\frac{1}{2}}  &0& \sqrt{x+\frac{1}{2}}
\end{pmatrix}\]
Using the expressions defined  in page (\pageref{e1}) we get
\begin{center}
$x_{0}=\frac{1-4}{2\sqrt{(4-1)^{2}+4 \times 2^{2}}}=\frac{-3}{10}$
\end{center}
\begin{center}
$\lambda_{1}=\frac{5}{2}+\frac{{(1-4)}^{2}+4 \times 2^{2}}{2\sqrt{(4-1)^{2}+4 \times 2^{2}}}=5$
\end{center}
\begin{center}
$\lambda_{2}=3$
\end{center}

\begin{center}
$\lambda_{3}=\frac{5}{2}+\frac{-{(1-4)}^{2}-4 \times 2^{2}}{2\sqrt{(4-1)^{2}+4 \times 2^{2}}}=0$
\end{center}
And the corresponding eigenvectors are

\[v_{1}=
\begin{pmatrix} 5\sqrt{5}\\ 0\\ 0 \end{pmatrix}, v_{2}=\begin{pmatrix}0\\ 3\\ 0 \end{pmatrix}, v_{3}=\begin{pmatrix} 0\\ 0\\ 0 \end{pmatrix} \]
 
\section{Conclusion}
In this paper, we gave another method for the computation of eigenvalues and eigenvectors of a real symetric matrix, and we well noticed the relationship between the two orthogonal matrices, i.e., these matrices allows us to calculate the same eigenvalues of a real symetric matrix but with two different values. Indeed, in the Jacobi method this value $\theta \in ] -\pi/4,\pi/4[$ and in this new method $x \in ]-\delta,1-\delta[$ such that $\delta \in \mathbb{R}$, so we can deduce that  there is a bijection between these two intervals.


\newpage

\begin{thebibliography}{160}
\bibitem{luca} L. Amodei and J P.  Dedieu,
{\em Analyse num\'erique matricielle}, Paris, Dunod, 2008.

\bibitem{fran} J. G. F. Francis,
{\em The QR transformation}, Part 1 and Part 2, Computer Journal, {\bf 4}, pp 265-271, pp 332-345, 1961, 1962.
\bibitem{jacobi} C. G. J. Jacobi  {\em $\ddot{U}$ber ein leichtes Verfahren, die in der Theorie der S$\ddot{a}$cularst\H{o}rungen vorkommenden Gleichungen numerisch aufzul\H{o}sen}, J. Reine Angew. Math. 30, pp 51-94, 1846.
\bibitem{frank} F. Jedrzejewski,
{\em Introduction aux m\'ethodes num\'eriques}, Deuxi\'eme \'edition, Springer.
\bibitem{kubla} V. N. Kublanovskaya,
{\em on some algorithms for the solution of the complete eigenvalue problem}, USSR. Compt. Math. Math. Phys, pp 637-657, 1961.
\bibitem{muntz} C. M$\ddot{u}$ntz,
{\em Solution directe de l'\'equation s\'eculaire et de quelques probl\`emes analogues transcendants}, Compte Rendu Acad, Paris, pp 43-46, 1913.
\bibitem{quart} A. Quarteroni, R. Sacco and  F. Saleri,
{\em M\'ethodes num\'eriques:  Algorithmes, Analyse et Applications}, Springer. 





 






\end{thebibliography}
\end{document}